\documentclass[11pt,letterpaper]{article}
\usepackage[T1]{fontenc}
\usepackage{lmodern,amsmath,amssymb,amsthm,microtype}
\usepackage[margin=1in]{geometry}
\usepackage[hidelinks,bookmarksnumbered]{hyperref}
\numberwithin{equation}{section}
\newtheorem{theorem}{Theorem}[section]
\newtheorem{proposition}[theorem]{Proposition}
\newtheorem{lemma}[theorem]{Lemma}
\newtheorem{corollary}[theorem]{Corollary}
\newcommand{\R}{\mathbb R}
\newcommand{\N}{\mathbb N}
\newcommand{\ii}{\mathrm i}
\newcommand{\e}{\mathrm e}
\newcommand{\dd}{\,\mathrm d}
\newcommand{\cH}{\mathcal H}

\newcommand{\st}{N_*}
\title{How many arithmetic modes recover the sign of a high derivative of a theta correlation?}
\author{Wei Sun\\[3pt]
\small Department of Mathematics and Statistics\\
\small Concordia University\\
\small Montreal, H3G 1M8, Canada\\
\small\href{mailto:wei.sun@concordia.ca}{\texttt{wei.sun@concordia.ca}}}
\date{}
\makeatletter
\renewcommand{\@maketitle}{%
  \begin{center}
    {\LARGE\@title\par}
    \vspace{0.8em}
    {\large\begin{tabular}{c}\@author\end{tabular}\par}
  \end{center}}
\makeatother
\hypersetup{pdftitle={How many arithmetic modes recover the sign of a high derivative of a theta correlation?},pdfauthor={Wei Sun}}
\begin{document}
\maketitle
\begin{abstract}
Let \(\cH(\tau,t)=2|\xi(\tfrac12+\tau+\ii t)|^2\), and let \(\cH_N\)
denote the truncation of its theta representation to pairs of positive
integers \((m,n)\) with \(mn\le N\). Fix a nonempty compact height interval
\(I\). For each fixed \(N\) and all sufficiently large \(p\), the
\(2p\)-th \(\tau\)-derivative at \(\tau=0\) is negative for \(\cH_N\) and
positive for \(\cH\) throughout \(I\). As \(p\to\infty\), the least
cutoff \(N\) for which this derivative of \(\cH_{N'}\) is positive on
\(I\) for every \(N'\ge N\) is asymptotic to
\(0.6044658090\ldots\,p/\log p\). The law is uniform over intervals contained in
\(|t|\le K\sqrt{p/\log p}\) for each fixed \(K>0\).
Comparison with the omitted pairs \(m=1,\ n>N\) gives a real center
whose ceiling is the exact cutoff outside a shrinking neighborhood of
the integers. Within that neighborhood at most two adjacent candidates
remain. One cutoff also makes a block of consecutive high derivatives
of \(-\partial_\tau\cH_N/\tau\) with respect to \(\tau^2\) strictly negative throughout
\(0\le\tau\le\tfrac12\), using the continuous extension at zero.
All results are unconditional.
\end{abstract}
\noindent\textit{Research collaboration.} This paper is based on research
that Wei Sun and ChatGPT (Alex), an AI system developed by OpenAI, conducted
jointly from 17 August through 25 November 2026. Their respective
contributions are described at the end of the paper.



\section{Introduction and main results}\label{sec:definition}

The Riemann xi function is the entire continuation of
\(\xi(z)=\tfrac12z(z-1)\pi^{-z/2}\Gamma(z/2)\zeta(z)\) from
\(\operatorname{Re}z>1\), with \(\xi(0)=\xi(1)=\tfrac12\).
For real \(\tau,t\), put \(\cH(\tau,t)=2|\xi(\tfrac12+\tau+\ii t)|^2\).
The parameter \(\tau\) measures horizontal displacement from the critical
line at height \(t\).

Sun used a modulus representation of \(\xi\) to obtain an analytic
inequality equivalent to the Riemann hypothesis~\cite{SunInequality}.
A representation in terms of the Jacobi theta function then gave a
family of polynomial inequalities whose joint validity is equivalent to
that hypothesis~\cite{SunModulus}. Here we ask how much of a theta
representation one must retain to recover the sign of a high horizontal
derivative of \(\cH\).

\subsection{Arithmetic truncation of the correlation}\label{sub:arithmetic}

Throughout, \(\N=\{1,2,\ldots\}\) and all logarithms are natural.
Set \(P(z)=4z^2-6z\). For \(m\in\N\) and \(x\in\R\), define
\begin{equation}\label{eq:phi}
 \phi_m(x)=\e^{x/4}P(\pi m^2\e^x)\e^{-\pi m^2\e^x},
 \qquad \phi(x)=\sum_{m\ge1}\phi_m(x).
\end{equation}
The complete correlation has the integral representation
\[
 \cH(\tau,t)=4\int_0^\infty\!\int_0^\infty
 \cos(tx)\cosh(\tau y)\phi(x+y)\phi(x-y)\dd x\dd y .
\]
Indeed, evenness of the complete kernel gives an integral over \(\R^2\)
with weight \(\e^{-\ii tx+\tau y}\). The substitution \(u=x+y,\ v=x-y\)
gives one half the product of two Mellin transforms of \(\phi\).
Section~\ref{sec:identities} evaluates these factors.

For \(N\in\N\), retaining all pairs whose product is at
most \(N\) gives the finite correlation
\begin{equation}\label{eq:truncation}
 \cH_N(\tau,t)=4\int_0^\infty\int_0^\infty
  \cos(tx)\cosh(\tau y)\sum_{mn\le N}\phi_m(x+y)\phi_n(x-y)\dd x\dd y .
\end{equation}
The \(k\)-th arithmetic mode groups all ordered pairs with \(mn=k\).
For \(x\ge0\), set
\[
 H_{\tau,k}(x)=\int_\R\cosh(\tau y)
             \sum_{m\mid k}\phi_m(x+y)\phi_{k/m}(x-y)\dd y,\qquad
 Q_{\tau,k}(x)=\partial_\tau H_{\tau,k}(x).
\]
As \(m\) runs through the positive divisors of \(k\), each ordered pair
occurs once. Interchanging \(m\) and \(k/m\) makes the divisor sum even
in \(y\), so its integral against \(\cosh(\tau y)\) is twice the
integral over \(y\ge0\). Extend each mode evenly in \(x\) and write
\(\widehat f(t)=\int_\R f(x)\e^{-\ii tx}\dd x\).
Then \(\cH_N=\sum_{k\le N}\widehat H_{\tau,k}\) and
\(\partial_\tau\cH_N=\sum_{k\le N}\widehat Q_{\tau,k}\).
The sum over all \(k\) gives \(\partial_\tau\cH=\widehat Q_\tau\), where
\(Q_\tau=\partial_\tau\sum_{k\ge1}H_{\tau,k}\). In~\eqref{eq:truncation}
the individual factors retain the signed arguments \(x+y\) and \(x-y\).
The identity \(\phi(x-y)=\phi(y-x)\) holds for the complete theta kernel.

\subsection{The stable cutoff}\label{sub:stable}

For an integer \(p\ge1\), write
\[
 M_{p,N}(t)=\partial_\tau^{2p}\cH_N(0,t),\qquad
 M_p(t)=\partial_\tau^{2p}\cH(0,t).
\]
Differentiation replaces \(\cosh(\tau y)\) by \(y^{2p}\) in the two
integrals at \(\tau=0\), so \(M_{p,N}\) and \(M_p\) are moments of the
finite and complete kernels. For fixed \(p\), \(M_{p,N}\to M_p\) as
\(N\to\infty\). For fixed \(N\) and sufficiently large \(p\),
\(M_{p,N}(t)<0<M_p(t)\) throughout any fixed compact height interval
(Proposition~\ref{pr:fixed}).
How large must \(N\) be for the positive sign to persist at every later
cutoff? For a nonempty compact set \(I\subset\R\), including a singleton, define
\begin{equation}\label{eq:stable}
 \st(p,I)=\min\{N\in\N:\ M_{p,N'}(t)>0\quad
                    \text{for every integer }N'\ge N\text{ and every }t\in I\}.
\end{equation}
We assign \(+\infty\) to the empty minimum. We use strict positivity
because the complete derivative is strictly positive throughout the
height range of our large-order result.

\begin{theorem}[Growth of the required cutoff]\label{th:main}
Let \(u_0\) be the positive solution of \(\e^{u_0}=2+u_0\) and put
\(h_0=\sqrt{u_0(2+u_0)}\). For every fixed nonempty compact interval
\(I\), as \(p\to\infty\),
\begin{equation}\label{eq:law}
 \st(p,I)\sim \frac{h_0}{\pi}\frac{p}{\log p},\qquad
      \frac{h_0}{\pi}=0.6044658090\ldots .
\end{equation}
The equivalent is uniform over all such \(I=I_p\) contained in
\([-K\sqrt{p/\log p},K\sqrt{p/\log p}]\) for each fixed \(K\).
\end{theorem}

The constant \(h_0/\pi\) comes from comparing the exponential rates of
the complete moment and the omitted row \(m=1,\ n>N\). In
\hyperref[sub:growthproof]{Subsection~\ref*{sub:growthproof}} we prove
Theorem~\ref{th:main} by testing one cutoff on each side of the point
where these rates agree. At the lower cutoff the signed omitted
moment exceeds the complete moment, so the truncated derivative is
negative. At the upper cutoff an absolute bound for the omitted terms
is smaller than \(M_p(t)\) throughout \(I\). This bound decreases as
the cutoff increases, giving positivity at every later index.

Which integer is the first stable cutoff for the derivative of order \(2p\)?
A finer comparison of the same two moments
gives a real center depending on the compact height set \(I\) through
\(T=\max_{t\in I}|t|\). Theorem~\ref{th:round} shows that the ceiling
of this center is the exact cutoff outside a shrinking neighborhood of
the integers. Within that neighborhood at most two adjacent candidates remain.

For \(\Xi(t)=\xi(\tfrac12+\ii t)\), the normalized complete moment
\(M_p(t)/(2p)!\) is the \(p\)-th generalized Laguerre expression for
\(\sqrt2\,\Xi\). Csordas~\cite{Csordas} represents these expressions as
Fourier transforms of weighted correlation kernels. Since
\(\cH(\tau,0)=2\xi(\tfrac12+\tau)^2\), \(M_p(0)\) is determined by the even
Taylor coefficients of \(\xi\) at \(1/2\), for which O'Sullivan~\cite{OSullivan}
and Romik~\cite{Romik} give asymptotic formulas. Haglund~\cite{Haglund}
studies the zeros of partial sums of the incomplete gamma expansion of
\(\Xi\) obtained by truncating a single theta kernel. The cutoff \(mn\le N\)
in~\eqref{eq:truncation} couples two theta factors. Comparing the omitted
moment with \(M_p(t)\) gives the sign of \(M_{p,N}(t)\) as \(p\) and \(N\)
grow together.

The remainder of this paper is organized as follows.
Section~\ref{sec:identities} derives the correlation and introduces the
two contact functions. The moment estimates in Sections~\ref{sec:full}
and~\ref{sec:tail} give the growth law (Theorem~\ref{th:main}) in
Subsection~\ref{sub:growthproof} and the rounding rule
(Theorem~\ref{th:round}) in Subsection~\ref{sub:rounding}.
Section~\ref{sec:strip} treats the full tilt interval, simultaneous orders
and the inverse problem. Section~\ref{sec:fixed} explains the negative sign
at a fixed cutoff, and Section~\ref{sec:numerics} illustrates the transition
numerically. The formal verification and research contributions follow.

\section{The complete correlation and finite contacts}\label{sec:identities}
Set \(g(x)=\e^{x/4}\sum_{m\ge1}\e^{-\pi m^2\e^x}\).
Differentiating each term yields the summands in~\eqref{eq:phi}, so
\(\phi=4g''-g/4\). Jacobi's transformation for \(\theta_3\) at zero
argument~\cite[(20.7.32)]{DLMF} gives \(g(-x)-g(x)=\sinh(x/4)\),
so \(\phi\) is even.
The Gaussian tails in~\eqref{eq:phi} give double-exponential decay
of \(\phi\) and its fixed derivatives as \(x\to+\infty\). Evenness gives
the corresponding bounds as \(x\to-\infty\). For \(\operatorname{Re}z>1\), integrating by parts twice
and evaluating the gamma integral gives
\[
 \int_{\R}\phi(x)\e^{(z-1/2)x/2}\dd x
 =z(z-1)\pi^{-z/2}\Gamma(z/2)\zeta(z)=2\xi(z).
\]
At the negative end, theta inversion gives \(g,g'=O(\e^{-x/4})\),
which makes the boundary terms vanish. The integral is entire in \(z\),
and the evenness of \(\phi\) gives \(\xi(1-z)=\xi(z)\).
The complete product \(\phi(x+y)\phi(x-y)\) is even in each variable.
Thus its integral over \(\R^2\) with weight \(\e^{-\ii tx+\tau y}\)
equals the quadrant integral defining \(\cH\), including the factor four.
With \(u=x+y,\ v=x-y\), the inverse Jacobian is \(1/2\) and the exponent
is \((\tau-\ii t)u/2+(-\tau-\ii t)v/2\). For \(z=\tfrac12+\tau+\ii t\)
the Mellin formula therefore gives
\[
 \cH(\tau,t)=\tfrac12\,[2\xi(\bar z)]\,[2\xi(1-z)]=2|\xi(z)|^2 .
\]
Reflecting the complete factor \(\phi(x-y)=\phi(y-x)\) gives
\begin{equation}\label{eq:complete}
 \cH(\tau,t)=2\int_0^\infty\cosh(\tau y)\kappa(y,t)\dd y,
 \qquad \kappa(y,t)=\int_{\R}\phi(y+x)\phi(y-x)\cos(tx)\dd x .
\end{equation}

For the individual terms of~\eqref{eq:truncation} retain \(x-y\) as
the second argument. Set \(Y=\e^y,\ v=\e^x\) and
\[
 A=\pi m^2Y,\quad D=\pi n^2/Y,\quad S=v(m^2Y+n^2/Y).
\]
Then \(\phi_m(x+y)\phi_n(x-y)=v^{1/2}P(vA)P(vD)\e^{-\pi S}\).
For every \(\eta>0\),
\begin{equation}\label{eq:summable}
 \sum_{m,n\ge1}\int_0^\infty\int_0^\infty
 v^{1/2}|P(vA)P(vD)|\e^{-\eta\pi S}\dd x\dd y<\infty.
\end{equation}
Indeed \(S\ge\e^{x+y}\ge v\) and \(S\ge2mn\). One factor of the
exponential absorbs the polynomials, a second is summable in \(m,n\), and
a third is integrable in \(x,y\), even with the weight \(y^q\e^{cy}\)
for fixed \(q\ge0\) and \(c\in\R\). This justifies the expansions,
termwise differentiation and analytic continuation in \(\tau\) on compact
sets. For each fixed \(p\), it also gives \(M_{p,N}\to M_p\) uniformly in \(t\in\R\).

\paragraph{The contact quantities.}
At positive tilt, a contact is a zero of \(\partial_\tau\cH_N(\tau,t)\).
To study its multiplicity at fixed \(N,t\), write \(s=\tau^2\) and define
the analytic functions
\[
 C_N(s,t)=-\frac{\partial_\tau\cH_N(\tau,t)}{\tau},\qquad
 L_N(s,t)=2\partial_s C_N(s,t)
\]
by continuity at zero. Define \(C\) and \(L\) in the same way from
\(\cH\). At a positive-tilt contact \(C_N=0\),
\(\partial_\tau^2\cH_N=-sL_N\), by differentiation of
\(\partial_\tau\cH_N=-\tau C_N\). At the boundary the Taylor series gives
\begin{equation}\label{eq:contact}
 \partial_s^{p-1}C_N(0,t)=-\frac{(p-1)!}{(2p-1)!}M_{p,N}(t).
\end{equation}
Thus \(C_N(0,t)=-M_{1,N}(t)\), \(L_N(0,t)=-M_{2,N}(t)/3\), and
the growth law measures the cutoff needed for every higher boundary
derivative in this hierarchy. With \(f=\sqrt2\,\Xi\), the functional
equation gives \(\cH(\tau,t)=f(t+\ii\tau)f(t-\ii\tau)\).
The coefficient of \(\tau^{2p}\) in this product is \(M_p(t)/(2p)!\).
For example \(M_1(t)=2(f'(t)^2-f(t)f''(t))\).

\section{The full moment}\label{sec:full}
Csordas, Norfolk and Varga~\cite{CNV} used moments of the theta kernel to prove
the Tur\'an inequalities for the normalized even Taylor coefficients of \(\xi\).
Here the moment is taken against the correlation density in~\eqref{eq:complete}.
Write \(q=2p\) and define \(b,R,a\) by
\begin{equation}\label{eq:scales}
 b\e^b=p/\pi,\qquad R=\e^b=\frac{p}{\pi b},\qquad
 a=2\pi R=\frac{q}{b}.
\end{equation}
In particular \(b=W(p/\pi)\), the positive branch of the Lambert function,
and \(b\sim\log p\). All three parameters tend to infinity. We use
\(g_b(u)=b\log(1+u/b)\) for \(u>-b\), so that
\begin{equation}\label{eq:linearize}
 (b+u)^q=b^q\e^{a g_b(u)},\qquad g_b(u)\le u .
\end{equation}
The letter \(K\) below denotes an arbitrary fixed bound for
\(|t|/\sqrt a\). All uniform errors may depend on \(K\). For \(G\ge0\),
we write \(F\ll G\) if \(|F|\le CG\) and \(F\gg G\) if \(F\ge cG\),
with positive constants uniform in the stated ranges.

\begin{lemma}[Complete density]\label{le:density}
Let \(Y=\e^y\). Then
\begin{equation}\label{eq:density}
 \sup_{t\in\R}\left|
     \frac{\kappa(y,t)}{16\pi^4Y^4\e^{-2\pi Y}}
        -\e^{-t^2/(4\pi Y)}\right|
     \longrightarrow 0 \qquad (Y\to\infty).
\end{equation}
For every \(y\ge0\) and real \(t\),
\( |\kappa(y,t)|\le C Y^4\e^{-2\pi Y}\).
\end{lemma}
\begin{proof}
On \(z\ge0\), summation of the Gaussian tail in~\eqref{eq:phi} gives
\(\phi(z)\ll\e^{9z/4}\exp(-\pi\e^z)\). As \(z\to\infty\) the ratio to
\(4\pi^2\e^{9z/4}\exp(-\pi\e^z)\) tends to one. On \(|x|\le y\),
the product in~\eqref{eq:complete} is consequently bounded by
\(CY^{9/2}\e^{-2\pi Y\cosh x}\). Put \(x=w/\sqrt{Y}\) and remove
\(16\pi^4Y^{9/2}\e^{-2\pi Y}\). At each fixed \(w\) the resulting
product tends to \(\e^{-\pi w^2}\), with an integrable majorant
\(C\e^{-\pi w^2}\). The \(L^1(\dd w)\) error therefore tends to zero.
It bounds the error in the cosine transform for every real \(t\), since
\(|\cos(tw/\sqrt Y)|\le1\). The Gaussian transform gives
\(\exp(-t^2/(4\pi Y))\). For \(|x|>y\), evenness gives the bound
\(C\e^{9|x|/2}\exp(-2\pi\e^{|x|}\cosh y)\) on the product.
Its integral is \(O(\exp(-cY^2))\) as \(Y\to\infty\), by the change
\(v=\e^{|x|}\). Enlarging the constant on bounded \(y\) proves both claims.
\end{proof}

\begin{proposition}\label{pr:full}
As \(p\to\infty\), uniformly for \(|t|\le K\sqrt a\),
\begin{equation}\label{eq:fullmoment}
 M_p(t)=32\pi^4 b^q R^4 \e^{-a}
 \sqrt{\frac{2\pi}{a(1+1/b)}}\,\e^{-t^2/(2a)}(1+o(1)).
\end{equation}
\end{proposition}
\begin{proof}
Use \(M_p=2\int_0^\infty y^q\kappa(y,t)\dd y\) and put \(y=b+u\).
The exponent after removing \(b^q\) is \(a(g_b(u)-\e^u)\), whose unique
maximum is \(-a\) at \(u=0\). Its second derivative there is \(-a(1+1/b)\).
On \(|u|\le1\) the normalized phase \(g_b(u)-\e^u\) has second
derivative at most \(-\e^{-1}\) and uniformly bounded third derivative
for large \(b\). For \(|u|>1\), use \(g_b(u)\le u\). The function \(u-\e^u\)
has a strict gap from \(-1\) on this set and
\(\int_{\R}\e^{4u}\exp(u-\e^u)\dd u<\infty\). Thus the integrated
contribution there, divided by \(b^qR^4\), is \(O(\e^{-(1+c)a})\) for
some \(c>0\), by Lemma~\ref{le:density}.
On \(|u|\le1\), put \(z=\sqrt a\,u\). The curvature bound gives
\(\exp\{a[g_b(u)-\e^u+1]\}\le\exp(-z^2/(2\e))\). For bounded \(z\),
Taylor expansion gives
\[
  a[g_b(u)-\e^u+1]= -\tfrac12(1+1/b)z^2 + O(|z|^3/\sqrt a).
\]
Here \(\e^{4u}=1+o(1)\), and~\eqref{eq:density} gives the normalized
density as \(\exp(-t^2/(4\pi R\e^u))+o(1)=\exp(-t^2/(2a))+o(1)\)
uniformly for \(|t|\le K\sqrt a\). Integrating the quadratic approximation
with the Gaussian majorant above and \(\dd u=\dd z/\sqrt a\) gives the
factor \(\sqrt{2\pi/[a(1+1/b)]}\). The height factor is at least
\(\e^{-K^2/2}\), so the error in~\eqref{eq:fullmoment} is relative. In particular the complete moments are positive there
for all sufficiently large \(p\).
\end{proof}

\section{The first omitted row}\label{sec:tail}
On the quadrant in~\eqref{eq:truncation}, the leading contribution from
\(mn>N\) comes from the row \(m=1,\ n=N+1,N+2,\ldots\). We compare
this row with the complete moment as the cutoff ratio \(h=N/R\) varies.
For \(h>0\) and \(u>-b\), put
\begin{align}
 Q(u,h)&=\tfrac12(\e^u+h^2\e^{-u}), &
 F_b(u,h)&=g_b(u)-Q(u,h),\label{eq:F}\\
 U_b(h)&=\mathop{\rm argmax}_{u>-b}F_b(u,h), &
 G_b(h)&=F_b(U_b(h),h).\label{eq:G}
\end{align}
The maximum exists and is unique. Indeed \(F_b\to-\infty\) at both ends
and \(\partial_u^2F_b=-Q-1/[b(1+u/b)^2]<0\). Set
\begin{equation}\label{eq:constants}
 U=U_b(h),\quad Q_b(h)=Q(U,h),\quad
 B_b(h)=Q_b(h)+\frac{1}{b(1+U/b)^2},\quad d_b(h)=2\pi h\e^{-U}.
\end{equation}
In the following discussion \(h=N/R\) stays in a fixed compact interval
\(J\subset(0,\infty)\). The maximizers stay in a common bounded interval
because \(F_b(0,h)\) is bounded below and \(u-Q(u,h)\to-\infty\) at both ends
uniformly in \(h\in J\). Consequently all the quantities in
\eqref{eq:constants} other than \(U\) are bounded above and bounded away
from zero. On bounded \(u\)-intervals the first three derivatives of
\(g_b\) are uniformly bounded for sufficiently large \(b\).

Write \(E_{p,N}(t)=M_p(t)-M_{p,N}(t)\) for the signed omitted moment.
Dropping the cosine and taking absolute values of the two polynomial
factors gives the envelope
\[
 A_{p,N}=4\sum_{mn>N}\int_0^\infty\!\int_0^\infty y^q
       v^{1/2}|P(vA)P(vD)|\e^{-\pi S}\dd x\dd y .
\]
Its summands are nonnegative, so for every \(n\ge N\) and real \(t\),
\begin{equation}\label{eq:major}
                |E_{p,n}(t)|\le A_{p,n}\le A_{p,N} .
\end{equation}

\begin{proposition}[Omitted contribution]\label{pr:tail}
Uniformly for \(N/R=h\in J\) and \(|t|\le K\sqrt a\), both
\(E_{p,N}(t)\) and \(A_{p,N}\) are asymptotic to
\begin{equation}\label{eq:tail}
 64\pi^4 b^q R^4 h^4\e^{aG_b(h)}a^{-3/2}
 \frac{\sqrt{2\pi/B_b(h)}}{Q_b(h)(\e^{d_b(h)}-1)} .
\end{equation}
\end{proposition}
\begin{proof}
For each pair set \(u=\log(m^2Y/R)\), \(j'=mn/R\). Then
\[
 y=b+u-2\log m,\qquad \pi S=a\e^x Q(u,j'),\qquad j'>h.
\]
The exponent per \(a\), after removal of \(b^q\), is
\(g_b(u-2\log m)-\e^x Q(u,j')\). Points with \(y=0\) have zero moment
weight and may be omitted.

\textit{Localization.} There exist \(L>0\) and \(c>0\), depending only
on \(J\), such that the contribution in absolute value from outside
\begin{equation}\label{eq:central}
                  m=1,\qquad |u|\le L,\qquad 0\le x\le1
\end{equation}
is \(O(b^q\e^{a(G_b(h)-c)})\). To see that summation preserves the gap,
reserve the factor \(\e^{-\eta\pi S}\) of~\eqref{eq:summable} with a
sufficiently small fixed \(\eta>0\). If \(|u|>L\), the remaining exponent
per \(a\) is bounded by \(u-(1-\eta)Q(u,h)\), which tends to \(-\infty\)
uniformly for \(h\in J,\ 0\le\eta\le1/2\). Choose \(L\) large enough
for a fixed gap there and to contain every maximizer in its interior.
If \(|u|\le L\) and \(x>1\), then \((1-\eta)\e^x>1+c_1\) with a fixed
\(c_1>0\), giving a gap since \(Q\ge\min J\). Finally if \(|u|\le L\)
and \(m\ge2\), the mean value theorem and \(g_b'(w)=(1+w/b)^{-1}\) give
\(g_b(u)-g_b(u-2\log m)\ge (\log4)/(1+L/b)\). In fact
\[
 g_b(u-2\log m)-(1-\eta)\e^x Q(u,j')
 \le F_b(u,h)-\frac{\log4}{1+L/b}+\eta Q(u,h),
\]
since \(\e^xQ(u,j')\ge Q(u,h)\). On this bounded \(u\)-interval,
choose \(\eta\) small enough that the last term costs at most half
of the gain, uniformly in \(h\in J\) and large \(b\). The remaining factors are integrable
and summable by~\eqref{eq:summable}.

\textit{Evaluation of the row.} In~\eqref{eq:central} write \(n=N+j\),
\(j\ge1\). Exactly,
\begin{equation}\label{eq:tailcost}
 \pi S=a\e^x Q(u,h)+\pi\e^{x-u}(2hj+j^2/R).
\end{equation}
For bounded \(j\), introduce \(z=\sqrt a\,(u-U)\), \(w=ax\). On each
bounded set of \((z,w,j)\) the full exponent, less \(aG_b(h)\), is
\begin{equation}\label{eq:localexp}
 -\tfrac12 B_b(h)z^2-Q_b(h)w-d_b(h)j+o(1).
\end{equation}
This follows by Taylor expansion at the stationary point, uniformly for
\(h\in J\). Both \(vA\) and \(vD\) tend uniformly to infinity there.
Since \(AD=\pi^2 n^2\) on this row, the factor
\(v^{1/2}P(vA)P(vD)/(16\pi^4R^4h^4)\) tends uniformly to one.
The same holds with absolute values. Also \(\cos(tx)=1+o(1)\), uniformly
for \(|t|\le K\sqrt a\).

On the whole central set, \(A,D\gg R\) and \(v\le\e\), so the
polynomial product is bounded in absolute value by \(C(AD)^2\).
Since \(AD=\pi^2(N+j)^2\) and \(N=Rh\), the normalized factor is
bounded by \(C(1+j/N)^4\le C(1+j)^4\). From \(\partial_u^2F_b\le-\min J\),
\(F_b(u,h)-G_b(h)\le-(\min J)(u-U)^2/2\). The term involving
\(\e^x-1\ge x\) in~\eqref{eq:tailcost} contributes at most
\(-a(\min J)x\) to the exponent. The cost of \(j\) is at least \(c_2j\)
for some \(c_2>0\). An integrable and summable majorant in the scaled
variables is thus \(C(1+j)^4\exp[-c_3(z^2+w+j)]\). Truncate these
variables to a fixed box, use the uniform local approximations there,
and then enlarge the box. The omitted tails are uniformly small by this
majorant. Integrating the exponential of the leading terms
in~\eqref{eq:localexp} over \(z\) and \(w\) gives the factors
\(\sqrt{2\pi/B_b(h)}\) and \(1/Q_b(h)\). Summation over \(j\ge1\) contributes
\((\e^{d_b(h)}-1)^{-1}\). The Jacobian is \(a^{-3/2}\).
Including the factor four in~\eqref{eq:truncation} proves~\eqref{eq:tail}
for both quantities. The prefactors are uniformly bounded away from zero
after the stated normalization, so the errors are relative.
\end{proof}

The integer sum in this calculation matters at the scale of one cutoff.
Its first term is \(j=1\), since the product \(mn=N\) is retained.
The complete moment has a Gaussian factor in \(t/\sqrt a\), whereas the
omitted row comes from \(x\) of order \(a^{-1}\). This difference explains
the height dependence of the transition.

\section{The growth law and its integer refinement}\label{sec:transition}
For \(h\in J\) define
\begin{equation}\label{eq:D}
 D_b(h)=\frac{2h^4\sqrt{1+1/b}}
   {Q_b(h)\sqrt{B_b(h)}(\e^{d_b(h)}-1)} .
\end{equation}
Propositions~\ref{pr:full} and~\ref{pr:tail} give the signed relative error
\begin{equation}\label{eq:ratio}
 \frac{E_{p,N}(t)}{M_p(t)}=
 \frac{D_b(h)}a\exp\left(a[G_b(h)+1]+\frac{t^2}{2a}\right)(1+o(1)).
\end{equation}
The same equivalent holds with \(A_{p,N}\) in the numerator. Both
relative errors are uniform for \(h\in J,\ |t|\le K\sqrt a\).

\subsection{Proof of Theorem~\ref*{th:main}}\label{sub:growthproof}

The limiting exponent for the omitted row is
\[
 G_\infty(h):=\max_{u\in\R}\{u-Q(u,h)\}
       =\log(1+\sqrt{1+h^2})-\sqrt{1+h^2}.
\]
Indeed the maximum occurs at \(\e^u=1+\sqrt{1+h^2}\). For \(h\) in
any fixed compact interval \(J\subset(0,\infty)\), this maximizer and
those in~\eqref{eq:G} lie in a common bounded \(u\)-interval, as shown
after~\eqref{eq:constants}. There \(g_b(u)\to u\) uniformly, so
\(G_b\to G_\infty\) uniformly on \(J\). The derivative
\(G_\infty'(h)=-h/(1+\sqrt{1+h^2})\) is negative. Since
\(\sqrt{1+h_0^2}=1+u_0\) and \(\log(2+u_0)=u_0\), we have
\(G_\infty(h_0)=-1\), matching the factor \(\e^{-a}\) in the complete
moment~\eqref{eq:fullmoment}.

For \(h=N/R\in J\), \(D_b(h)\) is bounded above and away from zero
by~\eqref{eq:constants} and~\eqref{eq:D}. Also \(t^2/(2a)\le K^2/2\)
for \(|t|\le K\sqrt a\). Both numerators and \(M_p(t)\) are positive
for large \(p\) by Propositions~\ref{pr:full} and~\ref{pr:tail}.
Taking logarithms in~\eqref{eq:ratio} and dividing by \(a\) gives the
uniform estimate
\[
 \frac1a\log\frac{Z}{M_p(t)}
       =1+G_\infty(N/R)+o(1),
 \qquad Z=E_{p,N}(t)\ \text{or}\ Z=A_{p,N}.
\]
Fix \(0<\delta<h_0\) and let \(N_\pm\) be the integer parts of
\(R(h_0\pm\delta)\). As \(N_\pm/R\to h_0\pm\delta\), the limiting
rate is positive at \(N_-\) and negative at \(N_+\). Hence
\(E_{p,N_-}(t)>M_p(t)>A_{p,N_+}\) throughout the height window for
all sufficiently large \(p\). Since \(M_{p,N_-}(t)<0\), every proposed
cutoff \(N\le N_-\) fails the sign test at index \(N_-\). For every
\(n\ge N_+\), the absolute bound~\eqref{eq:major} gives
\[
 M_{p,n}(t)\ge M_p(t)-A_{p,N_+}>0 .
\]
Thus \(N_-<\st(p,I_p)\le N_+\) for every nonempty interval \(I_p\)
in the window of the theorem, since \(a=2p/b\sim2p/\log p\). Every
fixed compact interval eventually lies in such a window.
Dividing by \(R\) and letting \(\delta\downarrow0\) gives
\(\st(p,I_p)/R\to h_0\) uniformly. Now \(R=p/(\pi b)\) and
\(b\sim\log p\) give~\eqref{eq:law}.\qed

\newpage
\subsection{Rounding the cutoff to an integer}\label{sub:rounding}

To distinguish adjacent cutoffs, we retain the finite-\(b\) exponent
and the factors multiplying the exponential in~\eqref{eq:ratio}.
Let \(u_b\) be the positive solution of the first equation below and
define \(h_b\) by the second.
\begin{equation}\label{eq:root}
\begin{aligned}
 \e^{u_b} &= 1+b\log(1+u_b/b)+(1+u_b/b)^{-1},\\
 h_b &= \sqrt{\e^{u_b}\{1+b\log(1+u_b/b)-(1+u_b/b)^{-1}\}} .
\end{aligned}
\end{equation}
The first equation has exactly one positive solution. Its left side minus
its right side is \(-1\) at zero, has derivative
\(\e^u-(1+u/b)^{-1}+b^{-1}(1+u/b)^{-2}>0\) for \(u\ge0\), and tends
to infinity. The expression under the square root is positive.
The equations in~\eqref{eq:root} give
\(\partial_uF_b(u_b,h_b)=0\) and \(F_b(u_b,h_b)=-1\). Moreover
\(u_b\to u_0\), the positive root of \(\e^u=2+u\), by local uniform
convergence and strict monotonicity across the limiting root. Thus \(h_b\to h_0\).
Set, for \(v\ge0\),
\begin{equation}\label{eq:center}
 \mathcal N_p(v)=R h_b+
 \frac{\log D_b(h_b)-\log a+v^2/(2a)}{d_b(h_b)} .
\end{equation}
All quantities on the right are functions of \(p,v\) only.
We write \(\lfloor x\rfloor\) for the greatest integer at most \(x\)
and \(\lceil x\rceil\) for the least integer at least \(x\).

\begin{theorem}[Integer transition]\label{th:round}
Fix \(K\ge0\). There exist \(p_0\in\N\) and a positive sequence \(\epsilon_p\to0\)
such that for every \(p\ge p_0\) and every nonempty compact
\(I_p\subset[-K\sqrt a,K\sqrt a]\), with \(T_p=\max_{t\in I_p}|t|\),
\begin{equation}\label{eq:round}
 \lfloor\mathcal N_p(T_p)-\epsilon_p\rfloor+1
 \ \le\ \st(p,I_p)\ \le\ \lceil\mathcal N_p(T_p)+\epsilon_p\rceil .
\end{equation}
This includes singleton heights. If the center stays a fixed positive
distance from the integers, the cutoff eventually equals
\(\lceil\mathcal N_p(T_p)\rceil\). Uniformly over the same sets,
\begin{equation}\label{eq:bounded}
        \st(p,I_p)=R h_b-\frac{\log a}{d_b(h_b)}+O_K(1).
\end{equation}
\end{theorem}
\begin{proof}
By stationarity, \(G_b'(h)=-h\e^{-U_b(h)}\). In a fixed neighborhood of
\(h_0\), the implicit function theorem, \(B_b\ge\min J>0\), and the
uniform compact bounds on \(U_b\) show that \(G_b''\) and
\((\log D_b)'\) are bounded independently of large \(b\). Take
\(N=R h_b+\Delta\) with \(\Delta=O(\log a)\). Since \(a/R=2\pi\),
\[
 a[G_b(N/R)+1]=-d_b(h_b)\Delta+O((\log a)^2/a).
\]
On the same range,
\(\log D_b(N/R)=\log D_b(h_b)+O((\log a)/a)\).
Equation~\eqref{eq:ratio} then gives, for either \(Z=E_{p,N}(t)\) or
\(Z=A_{p,N}\),
\begin{equation}\label{eq:profile}
 \log\frac{Z}{M_p(t)}=-d_b(h_b)\,[N-\mathcal N_p(|t|)]+o(1)
\end{equation}
uniformly for \(|N-Rh_b|\le C\log a\) and \(|t|\le K\sqrt a\) with
\(C,K\) fixed. The \(o(1)\) includes the integral error in
Proposition~\ref{pr:full} and the integration and summation errors in
Proposition~\ref{pr:tail}. The latter also gives \(E_{p,N}(t)>0\) here for
large \(p\). Since \(d_b(h_b)\) is bounded away from zero, choose
\(\epsilon_p\downarrow0\) slowly enough that
\(d_b(h_b)\epsilon_p\) exceeds twice the uniform absolute error in
\eqref{eq:profile} for \(|\Delta|\le C_K\log a\), taking \(C_K\) large
enough to include the two integers used below.

Set \(N_-=\lfloor\mathcal N_p(T_p)-\epsilon_p\rfloor\). At a height
with absolute value \(T_p\), \(E_{p,N_-}>M_p\), hence \(M_{p,N_-}<0\).
Set \(N_+=\lceil\mathcal N_p(T_p)+\epsilon_p\rceil\). Since the center
increases with \(|t|\), \(A_{p,N_+}<M_p(t)\) for every \(t\in I_p\).
Equation~\eqref{eq:major} gives \(M_{p,n}(t)>0\) for every \(n\ge N_+\)
and every such \(t\). This proves~\eqref{eq:round}. If the center is farther
than \(\epsilon_p\) from the nearest integer, both bounds equal its ceiling.
If \(|\mathcal N_p(T_p)-j|\le\epsilon_p<1/2\) for an integer \(j\), they
allow only \(j\) and \(j+1\). Finally, \(\log D_b(h_b)\) and \(T_p^2/(2a)\)
are bounded, giving~\eqref{eq:bounded}.
\end{proof}

\paragraph{The first correction to \(h_0\).}
The limiting values are \(u_0=1.1461932206\ldots\) and
\(h_0=1.8989853449\ldots\). To find the correction at finite \(b\),
write \(\beta=1/b\). The equation
for \(u_b\) is smooth at \(\beta=0\) and its \(u\)-derivative at \(u_0\)
is \(\e^{u_0}-1>0\). Hence \(u_b,h_b\) are smooth in \(\beta\) there.
Differentiate \(F_b(u_b,h_b)=-1\) at \(\beta=0\). The
\(u\)-derivative vanishes, while \(\partial_\beta g_b(u)|_{\beta=0}=-u^2/2\)
and \(\partial_hF_b=-h\e^{-u}\). Using \(h_0^2=u_0\e^{u_0}\) gives
\[
 h_b=h_0-\frac{h_0u_0}{2b}+O(b^{-2}).
\]
Thus replacing \(h_b\) by \(h_0\) in~\eqref{eq:bounded} would shift
the center by order \(R/b\), larger than \(\log a\).

\section{A common cutoff for the two contact functions}\label{sec:strip}
So far the recovered signs have been those of derivatives at zero tilt.
Ordinary contact tests use \(C_N(s,t)\) in the interior, and tests of
nearly degenerate contacts also use \(L_N(s,t)\) and higher derivatives
at the same point. We next allow the evaluation point to range over the
entire interval \(0\le s\le1/4\) and require a common cutoff.

For \(p\ge1\) let
\begin{equation}\label{eq:Z}
 Z_{p,N}(s,t)=-\frac{(2p-1)!}{(p-1)!}\,\partial_s^{p-1}C_N(s,t).
\end{equation}
Define \(Z_p\) by the same formula with \(C\) in place of \(C_N\).
At \(s=0\), identity~\eqref{eq:contact} gives \(Z_{p,N}(0,t)=M_{p,N}(t)\)
and \(Z_p(0,t)=M_p(t)\). At arbitrary \(s\), \(Z_{1,N}=-C_N\)
and \(Z_{2,N}=-3L_N\). Define the positive factor
\[
 W_p(s,y)=\sum_{j=0}^\infty\frac{(s y^2/4)^j}{j!(p+1/2)_j},\qquad
 (v)_j=v(v+1)\cdots(v+j-1),\quad (v)_0=1 .
\]
In each integral for \(M_{p,N}\), \(M_p\) or the signed difference, multiply
the integrand by \(W_p(s,y)\) to obtain the respective formula for \(Z\).
Indeed the coefficient of \(s^j\) after differentiating
\(y\sinh(\sqrt s\,y)/\sqrt s\) exactly \(p-1\) times and applying the
positive normalization in~\eqref{eq:Z} is
\[
 \frac{(2p-1)!}{(p-1)!}\frac{(p-1+j)!}{j!(2p-1+2j)!}\,y^{2p+2j}
 =\frac{y^{2p+2j}}{4^j j!(p+1/2)_j} .
\]
For \(0\le s\le1/4\), termwise comparison gives
\begin{equation}\label{eq:Wweight}
 1\le W_p(s,y)\le\min\left\{\exp\!\left(\frac{s y^2}{4p+2}\right),
                                       \cosh(\sqrt s\,y)\right\}.
\end{equation}
The two majorants use \((p+1/2)_j\ge(p+1/2)^j\) and
\((p+1/2)_j\ge(1/2)_j\), respectively.

\begin{theorem}[Uniformity in the tilt]\label{th:strip}
Formulas~\eqref{eq:fullmoment}, \eqref{eq:tail} and~\eqref{eq:ratio} hold
with \(M_p,E_{p,N},A_{p,N}\) replaced by \(Z_p, Z_p-Z_{p,N}\) and the
absolute envelope with weight \(W_p(s,y)\). All relative errors are
uniform for \(0\le s\le1/4\) together with the stated height and cutoff
ranges. The rounding bound~\eqref{eq:round} holds, with a possibly
enlarged error sequence \(\epsilon_p\to0\), for the least cutoff
\(N\) such that \(Z_{p,n}(s,t)>0\) for every \(n\ge N\),
\(s\in[0,1/4]\) and \(t\in I_p\).
\end{theorem}
\begin{proof}
On \(y=b+O(1)\), the first bound in~\eqref{eq:Wweight} gives
\(W_p=1+O(b^2/p)=1+o(1)\). This covers the complete saddle and the
entire central set~\eqref{eq:central} in the dominant row. On \(|u|\le1\)
the integral of \(y^q|\kappa(y,t)|\) is
\(O(b^qR^4\e^{-a}/\sqrt a)\) by Lemma~\ref{le:density} and the curvature
bound in Proposition~\ref{pr:full}. Thus multiplying the signed integrand
by \(W_p\) changes the complete central contribution by
\(O_K(b^2/p)\) times the leading term in~\eqref{eq:fullmoment}.
In the complete outer estimate \(|u|>1\), the second bound contributes at most
\(R^{1/2}\e^{u/2}\). The integral of
\(\e^{(9/2)u}\exp(u-\e^u)\) is finite, and the polynomial factor
\(R^{1/2}\) is absorbed by the existing exponential gap. Outside the
central set of the arithmetic tail the same bound costs at most
\(\e^{y/2}\), already permitted in~\eqref{eq:summable}. This proves the
two asymptotics with the weight and hence their ratio. The weighted
absolute envelope decreases in \(N\) for every \(s\). The upper argument
of Theorem~\ref{th:round} is therefore uniform in \(s\), while the lower
obstruction can already be tested at \(s=0\).
\end{proof}

\begin{corollary}[Simultaneous derivative orders]\label{co:joint}
For integers \(q_1\ge p\), require all the signs \(Z_{j,n}(s,t)>0\)
for \(p\le j\le q_1\), every cutoff \(n\ge N\), every \(0\le s\le1/4\),
and all \(t\) in a nonempty compact \(I_p\subset[-K\sqrt{a},K\sqrt{a}]\).
The least such \(N\)
satisfies~\eqref{eq:round} with center \(\mathcal N_{q_1}(T_p)\) and with
a vanishing error uniform over \(q_1\ge p\). In particular, \(q_1=p+1\)
recovers the signs of \(\partial_s^{p-1}C_n\) and \(\partial_s^{p-1}L_n\)
on the same strip. Both derivatives are then strictly negative.
\end{corollary}
\begin{proof}
Write \(b=b_p\), \(a_b=2\pi\e^b\), and \(r_b=d_b(h_b)\).
Then \(r_b\to r_*=2\pi h_0\e^{-u_0}>0\), and \(D_b(h_b)\)
has a positive limit by~\eqref{eq:constants} and its defining formula.
We compare two orders at a fixed height \(v\), with
\(v^2/a_b\le K^2\). Increasing the order increases \(b\), since
\(\pi b\e^b=p\).

For fixed \(u\ge0\), the quantity \(b\log(1+u/b)\) is nondecreasing
in \(b>0\). Its derivative is
\(\log(1+u/b)-1+(1+u/b)^{-1}\ge0\).
The term \((1+u/b)^{-1}\) also increases with \(b\).
In the scalar equation for \(u_b\), the difference between the two
sides is strictly increasing in \(u\), so \(u_c\ge u_b\) for \(c\ge b\).
Also \(h_c\ge h_b\). Indeed, if \(h_c<h_b\), evaluating the
phase \(F_c(\cdot,h_c)\) at \(u_b\) would give a value larger than
\(F_b(u_b,h_b)=-1\), contradicting its maximum of \(-1\).

Stationarity gives
\[
 r_b^2=(2\pi)^2\{1-2\e^{-u_b}(1+u_b/b)^{-1}\}.
\]
The function \(\e^{-u}(1+u/c)^{-1}\) decreases in \(u\).
At fixed \(u=u_b\le2\), the increase of \((1+u/b)^{-1}\)
when \(b\) is replaced by \(c\) equals
\(u(c-b)/[(b+u)(c+u)]\le2(c-b)\) for \(b\ge1\).
Consequently \(r_b^2-r_c^2\le4(2\pi)^2(c-b)\). Since the
spacings are bounded away from zero, there is a constant \(M\) with
\(r_c^{-1}-r_b^{-1}\le M(c-b)\) at all sufficiently large orders.

In~\eqref{eq:center} replace the bounded correction
\[
 (\log D_b(h_b)-\log(2\pi))/r_b+v^2/(2a_b r_b)
\]
by its first term's limit plus \(v^2/(2a_b r_*)\).
The error tends to zero uniformly for \(v^2/a_b\le K^2\).
The resulting center is nondecreasing with the order on this
window for all sufficiently large \(b\). To see this, the gain
in \(\e^b h_b\) is at least \(\e^b h_b(c-b)\).
The increase of \(c/r_c-b/r_b\) is bounded above by
\((bM+1/r_c)(c-b)\), and the decrease in the height term is at
most \(K^2(c-b)/(2r_*)\), using \(1-\e^{-(c-b)}\le c-b\).
The first bound dominates these two costs as \(b\to\infty\).

Thus every earlier center in a block is at most the terminal center
plus a uniform error tending to zero with the starting order.
The original height set belongs to every later window.
Apply Theorem~\ref{th:strip} at each order. The upper bound therefore
holds at every order in the block, and the terminal order supplies the
lower bound. Taking a common vanishing margin gives~\eqref{eq:round}.
At adjacent orders, \(L_n=2\partial_s C_n\) and the positive factorial
factors defining \(Z_{j,n}\) give the two strict signs.
\end{proof}

One consequence for contact geometry is an upper bound that counts
multiplicity. For all cutoffs beyond the sign threshold in
Theorem~\ref{th:strip} and all heights in its range, \(s\mapsto C_n(s,t)\)
has at most \(p-1\) zeros on \([0,1/4]\). For \(p\ge2\), \(L_n(s,t)\)
has at most \(p-2\) zeros there. Indeed \(\partial_s^{p-1}C_n<0\) and
\(\partial_s^{p-2}L_n=2\partial_s^{p-1}C_n\). Repeated Rolle's theorem
counts a zero of multiplicity \(m\) as \(m-1\) zeros at that point after
one differentiation, in addition to the zeros between distinct roots.
Both endpoints are included, with multiplicities measured in \(s\).
For a zero of \(C_n(\cdot,t)\) at \(s>0\), the contact of
\(\partial_\tau\cH_n\) at \(\tau=\sqrt s\) has the same multiplicity.
At zero tilt, if \(C_n(\cdot,t)\) has order
\(r\ge0\) at \(s=0\), the identity
\(\partial_\tau\cH_n=-\tau C_n(\tau^2,t)\) gives order \(2r+1\) in \(\tau\).

\begin{corollary}[Derivative orders available at a given cutoff]\label{co:inverse}
Fix a nonempty compact height set \(I\). Choose \(p_0\) large enough
that the cutoff of Theorem~\ref{th:strip} is finite at every order
\(j\ge p_0\) on \(I\). For all sufficiently large \(N\), let \(P(N)\)
be the largest integer \(p\ge p_0\) for which
\(Z_{j,n}(s,t)>0\) simultaneously for \(p_0\le j\le p\), \(n\ge N\),
\(0\le s\le1/4\), and \(t\in I\). Then, as \(N\to\infty\),
\[
                    P(N)\sim \frac{\pi}{h_0}\,N\log N .
\]
\end{corollary}
\begin{proof}
Write \(c=h_0/\pi\). By Theorem~\ref{th:strip}, the individual cutoffs
at order \(j\) are asymptotic to \(cj/\log j\). Their maximum from
\(p_0\) through \(p\) has the same equivalent, since \(j/\log j\) is
eventually increasing and any fixed prefix has bounded maximum. For
\(0<\delta<1\), take \(p_\pm=\lfloor(1\pm\delta)N\log N/c\rfloor\).
The respective maxima are asymptotic to \((1\pm\delta)N\).
For large \(N\), the orders through \(p_-\) meet the sign condition,
and the orders through \(p_+\) include one that fails it. The admissible
orders therefore have a maximum \(P(N)\) with \(p_-\le P(N)<p_+\).
Letting \(\delta\downarrow0\) gives the conclusion.
\end{proof}

\section{Fixed cutoffs and the scope of the law}\label{sec:fixed}
The elementary obstruction that motivated the growing-cutoff problem can
also be seen directly in the lattice integral. For fixed \(N\) let
\[
 k_N(y,t)=4\sum_{mn\le N}\int_0^\infty \cos(tx)\phi_m(x+y)\phi_n(x-y)\dd x,
 \qquad M_{p,N}(t)=\int_0^\infty y^q k_N(y,t)\dd y .
\]
\begin{proposition}[A fixed truncation has a negative outer density]\label{pr:fixed}
Uniformly for \(t\) in any fixed compact interval, as \(Y=\e^y\to\infty\),
\begin{equation}\label{eq:negative}
 k_N(y,t)=-96\pi^2\left(\sum_{n=1}^N n^2\right)\e^{-\pi Y}(1+o(1)).
\end{equation}
Consequently \(M_{p,N}(t)<0<M_p(t)\) uniformly on the interval for all
sufficiently large \(p\).
\end{proposition}
\begin{proof}
For one fixed pair set \(w=A(v-1)\). Then \(\dd x=\dd w/(A+w)\) and
\(\pi S=A+w+D(1+w/A)\). Write \(\lambda=1/A\), \(c=AD=\pi^2m^2n^2\),
so \(v=1+\lambda w\).
Since \(P(z)=z(4z-6)\), the polynomial factors and the Jacobian give
\[
 \frac{P(A+w)P(Dv)}{A+w}
   =(4v-6\lambda)(cv)(4c\lambda v-6).
\]
After multiplication by \(v^{1/2}\e^{-c\lambda v}\), this tends to \(-24c\)
for each fixed \(w\ge0\) as \(Y\to\infty\). After factoring out
\(\e^{-A}\), the remaining factors are \(\e^{-w}\) and
\(\cos(t\log v)\), the latter tending to \(1\) uniformly for bounded
\(t\). When \(0<\lambda\le1\), one has \(1\le v\le1+w\). The absolute value
of the rescaled integrand is at most
\[
   10c(4c+6)(1+w)^4\e^{-w},
\]
an integrable bound independent of \(\lambda\) and \(t\).
Dominated convergence along every sequence of bounded heights gives
the uniform pair limit. Finite summation gives~\eqref{eq:negative}, since \(m=1\) has
the slowest exponential. For uniform negativity of the moments choose
\(y_0>0\) beyond which the density is uniformly negative. The absolute
contribution below \(y_0\) is \(O(y_0^q)\). A fixed subinterval of
\([y_0+1,y_0+2]\) contributes at most \(-c(y_0+1)^q\). The rest beyond
\(y_0\) has the same sign. Proposition~\ref{pr:full} gives the complete sign.
\end{proof}

There is also a short proof of the complete sign on a fixed bounded
height interval directly from Lemma~\ref{le:density}. For large \(y\)
the supremum in~\eqref{eq:density} is less than \(1/4\), and the Gaussian
exceeds \(1/2\) throughout the interval.
Thus \(\kappa\) has a uniform positive lower bound on a fixed interval
\([y_0+1,y_0+2]\) and is positive for all \(y\ge y_0\). Comparing moments
on the two sides of \(y_0\) as above proves the claim.
For every fixed positive integer \(B\), the negative sign at \(N=B\)
also forces \(\st(p,I)>B\) for large \(p\), uniformly over nonempty
height sets in the bounded interval. Every proposed cutoff at or below
\(B\) must pass the test at \(B\) itself. This argument includes the
value \(+\infty\) in~\eqref{eq:stable}.

Ki~\cite{Ki} shows that the high even derivatives of \(\Xi\), after
normalization and rescaling of the argument, converge locally uniformly
to the cosine. For the correlation moment \(M_p(t)\),
formula~\eqref{eq:fullmoment} gives the Gaussian factor
\(\exp(-t^2/(2a))\) on the height window of Theorem~\ref{th:round}.
When \(|t|/\sqrt a\to\infty\), this Gaussian tends to zero. Extending the
argument to those heights would require control of the Fourier error
relative to the Gaussian in the complete saddle region.

\section{Numerical illustration}\label{sec:numerics}
For each row of Table~\ref{tab:transition}, take \(I=\{t\}\) and write
\(k=\lfloor\mathcal N_p(|t|)\rfloor\). Rational interval arithmetic
certifies each center to six decimal places and determines its floor.
The last two columns give quadrature estimates of \(E_{p,N}(t)/M_p(t)\)
at \(N=k\) and \(N=k+1\), without rigorous error enclosures.
When \(M_p(t)>0\), a ratio above
one gives \(M_{p,N}(t)<0\), and a ratio below one gives \(M_{p,N}(t)>0\).
For example, at \(p=256,\ t=0\) the computed ratio drops from
\(4.441770\) at \(N=40\) to \(0.087965\) at \(N=41\).

\begin{table}[ht]
\centering
\caption{Transition centers, floors, and signed moment ratios at
adjacent cutoffs.}\label{tab:transition}
\smallskip
\begin{tabular}{rcrrrr}
\hline
 & & \multicolumn{2}{c}{Interval arithmetic}
   & \multicolumn{2}{c}{Numerical quadrature}\\
\cline{3-4}\cline{5-6}
\(p\)&\(t/\sqrt a\)&\(\mathcal N_p(|t|)\)&\(k\)&
\(E_{p,k}/M_p\)&\(E_{p,k+1}/M_p\)\\
\hline
256  &0& 40.400401& 40& 4.441770&0.087965\\
1024 &0&127.191111&127& 2.051327&0.039814\\
8192 &0&752.874624&752&30.618936&0.606084\\
8192 &1&753.001829&753& 0.998559&0.019717\\
\hline
\end{tabular}
\end{table}

\begin{samepage}
In the last row, \(p=8192\) and \(t=\sqrt a\). The certified center lies
just above \(753\), so its ceiling is \(754\). At \(N=753\) the computed
signed ratio is \(0.998559\).
To certify \(\st(p,\{t\})=753\), it would suffice to prove
\(E_{p,752}(t)>M_p(t)>A_{p,753}\).
The first inequality would force a negative derivative at 752.
The second, by~\eqref{eq:major}, would give a positive derivative
at every cutoff \(n\ge753\).\par
\end{samepage}

For the complete integral we used Gauss--Legendre quadrature in
\(\sqrt{a}(y-b)\) and Gauss--Hermite quadrature in \(\sqrt{\pi Y}\,x\)
after reflecting each full theta factor. For the omitted integral we
used Gauss--Legendre quadrature in \(\sqrt{a}(\log(m^2Y/R)-U_b(N/R))\) and
Gauss--Laguerre quadrature in \((A+D)(\e^x-1)\). Increasing the outer
quadrature orders from 110 and 100 to 155 and 145, respectively,
and increasing the omitted sum from four rows with 35 terms each to
five rows with 48 terms each preserved the displayed digits of all
eight ratios.

For the four scalar centers, we bracketed \(b\) and the positive
root \(u_b\) in~\eqref{eq:root} by bisection, using the exponential
series and the identity
\(\log y=2\sum_{j\ge0}z^{2j+1}/(2j+1)\), where \(z=(y-1)/(y+1)\).
After rescaling \(y\) to \([1,2]\), the omitted terms have a geometric
majorant. The resulting intervals each lie within one rounding interval
for the printed decimal and between two consecutive integers.

\section*{Methods and scope of the formal verification}\label{sec:formal}
\addcontentsline{toc}{section}{Methods and scope of the formal verification}
Wei Sun used OpenAI's ChatGPT (Alex, GPT-6-based Codex) to develop the
moment estimates in Propositions~\ref{pr:full} and~\ref{pr:tail} and the
contact-derivative estimates in Section~\ref{sec:strip}. Alex wrote the
Lean proof scripts for this paper.
Lean~4.24.0 checked unconditional proofs of all nine numbered results
for the functions defined here. The proof of Theorem~\ref{th:main} bounds
the least stable cutoff~\eqref{eq:stable} uniformly over the stated height
sets and proves positivity at every later index. Exact bounds give
\(0.6044658090<h_0/\pi<0.6044658091\).

For Proposition~\ref{pr:tail} the verification identifies the signed
omitted sum with the difference of the two original derivatives. On the
first omitted row it verifies the change of variables, including the
Jacobian and the polynomial factors. A Gaussian in \(z\), exponential
decay in \(w\), and a summable bound in the integer offset \(j\) justify
integration and summation of the limit. An exponential bound controls
the remaining lattice points. The proof treats both the signed sum and
its absolute majorant.

The formal proof of Theorem~\ref{th:round} estimates the quadratic
remainder in the maximum directly by comparing its values at two cutoff
ratios. It proves~\eqref{eq:profile}, the common vanishing rounding
margin~\eqref{eq:round}, the exact ceiling away from the integers, and
the uniform bounded remainder~\eqref{eq:bounded}. For
Theorem~\ref{th:strip}, the proof identifies the normalized derivatives
with weighted integrals, including at zero tilt, and proves the same
rounding bounds for one cutoff on the entire tilt interval.
For Corollary~\ref{co:joint}, comparison of the finite balance heights
bounds the earlier centers by the last center plus a uniform vanishing
error. This gives the stated integer bounds for arbitrary blocks and
both strict signs at adjacent orders. For Corollary~\ref{co:inverse},
inverting the cutoffs on a fixed compact set gives the equivalent for
the largest admissible order. The zero count following Corollary~\ref{co:joint} is verified by
Rolle's theorem on ordered lists with repeated nodes. Local factorization
also verifies the contact orders under \(s=\tau^2\), including \(2r+1\)
at zero tilt.

Internal AI reviews clarified the complete-moment estimate and the bound
for every later cutoff, and independently reproduced the interval
calculations and quadrature estimates in Table~\ref{tab:transition}
using the methods of Section~\ref{sec:numerics}.

\section*{Contributions and use of generative AI}
\addcontentsline{toc}{section}{Contributions and use of generative AI}
\begin{samepage}
Wei Sun's work on modulus representations of \(\xi\)
\cite{SunInequality,SunModulus} motivated the joint research that he and
OpenAI's ChatGPT (Alex) carried out from 17 August through 25 November
2026. Sun directed the research and formulated this paper's question of
how the cutoff must grow with derivative order.
Alex developed the comparison between the complete moment and the
omitted lattice sum and checked constants and steps in the proofs.
Their joint work also produced Lean proofs of all nine numbered results
(\hyperref[sec:formal]{see the preceding section}).
\end{samepage}

\subsection*{Declaration of generative AI use in manuscript preparation}
\addcontentsline{toc}{subsection}{Declaration of generative AI use in manuscript preparation}
During manuscript preparation and revision (25--29 November 2026),
OpenAI's ChatGPT (Alex) drafted the exposition. Wei Sun reviewed the
manuscript, revised it jointly with Alex, and takes full responsibility
for the mathematics and the text.

\end{document}